\documentclass[12pt]{amsart}

\usepackage{amsmath,amsthm,amsfonts,amssymb}
\usepackage[all]{xy}

\numberwithin{equation}{section}

\theoremstyle{definition}

\DeclareMathOperator{\ev}{ev}

\DeclareMathOperator{\coh}{H}
\DeclareMathOperator{\Ext}{Ext}
\DeclareMathOperator{\Hom}{Hom}

\newcommand{\ot}{\otimes}
\newcommand{\V}{\mathcal{V}}

\newcommand{\DOT}{\setlength{\unitlength}{1pt}\begin{picture}(2.5,2)
                  (1,1)\put(2,3.5){\circle*{2}}\end{picture}}

\begin{document}

\title[Erratum]
{Erratum to ``Support varieties and representation type of small quantum groups''}

\author{J\"org Feldvoss}
\address{Department of Mathematics and Statistics,
University of South Alabama, Mobile, AL 36688--0002, USA}
\email{jfeldvoss@southalabama.edu}

\author{Sarah Witherspoon}
\address{Department of Mathematics, Texas A\&M University,
College Station, TX 77843--3368, USA}
\email{sjw@math.tamu.edu}

\date{August 22, 2013}

\maketitle


\begin{abstract}
Some of the general results in the paper require an additional hypothesis,
such as quasitriangularity.
Applications 
to specific types of Hopf algebras are correct, as some
of these are quasitriangular, and for those that are not,
the Hochschild support variety theory may be applied instead.
\end{abstract}


\vspace{.2cm}

Let $A$ be a finite dimensional Hopf algebra over a field $k$. 
The vector space $\coh^{*}(A,k):=\Ext^{*}_A(k,k)$ 
is an associative, graded commutative $k$-algebra under the cup
product, or equivalently under Yoneda composition. 
If $M$ and $N$ are finitely generated 
left $A$-modules, then $\coh^{*}(A,k)$ acts on $\Ext^{*}_A
(M,N)$ via the cup product, or equivalently by $-\ot N$ followed by Yoneda composition.

Let $\overline{S}$ be the composition inverse of the antipode $S$.
The isomorphism towards the top of p.\ 1350 in \cite{FW},
$\Ext^{\DOT}_A(M,N)\cong \Ext^{\DOT}_A(k, M^{*}\ot N)$,
assumes the following $A$-module structure on $M^* = \Hom_k(M,k)$,
which was not stated explicitly in the paper: 
$$
   (a\cdot f)(m) = \sum a_2 f(\overline{S}(a_1) m)
$$
for all $a\in A$, $f\in\Hom_k(M,k)$, and $m\in M$.

Under some finiteness assumptions as in \cite{FW}, 
we recall the definition of support variety: 
Let $M$ be a finitely generated left $A$-module. 
Let $I_A(M)$ be the annihilator of the action of
$\coh^{\ev}(A,k)$ on $\Ext^{\DOT}_A(M,M)$, 
and let $\V_A(M)$ denote the maximal ideal spectrum of the finitely generated
commutative $k$-algebra $\coh^{\ev}(A,k)/I_A(M)$.

In part (5) of \cite[Proposition 2.4]{FW}, it is stated that
\[
    \V_A(M\ot N) \subseteq \V_A(M)\cap \V_A(N).  \tag{*}\label{tp}
\]
Our proof of this statement requires an additional hypothesis
such as  quasitriangularity of $A$,
as we will explain. 
However, (\ref{tp}) holds  under other hypotheses,
for example, whenever $\V_A(M) = \V_A(M^*)$ for all finitely generated $A$-modules $M$,
as follows from the fact that dualization reverses the order of tensor products.
For a counterexample to the general statement (\ref{tp}),
see \cite{BW}. 

In the proof of (\ref{tp}) given in \cite{FW}, 
the action of $\coh^{ev}(A,k)$ on $\Ext^{\DOT}_A(M\ot N, M\ot N)$ does
indeed factor through its action on $\Ext^{\DOT}_A(M,M)$ as stated,
since we may first apply $-\ot M$, then $-\ot N$, then Yoneda composition. 
If $A$ is quasitriangular, then $M\ot N\cong N\ot M$, so  the action
factors through that on $\Ext^{\DOT}_A(N,N)$ as well, and consequently 
(\ref{tp}) holds. 
However, it does not necessarily factor through its action on $\Ext^{\DOT}_A(N,N)$
in general:
Under the isomorphism on $\Ext$ that is used in the proof,
$\Ext^{\DOT}_A(M\ot N, M\ot N) \cong \Ext^{\DOT}_A(N, M^{*}\ot M\ot N)$, 
the actions  of $\coh^{ev}(A,k)$ do not always coincide.
The order of the tensor product affects the action. 
Thus we obtain
$\V_A(M\ot N) \subseteq \V_A(M)$ in general,
but not necessarily (\ref{tp}).

Theorem 2.5,  Corollary 2.6,
Theorem 3.1, and Corollary 3.2 of \cite{FW} 
require an additional hypothesis such as quasitriangularity of $A$, 
since they use \cite[Proposition~2.4(5)]{FW}. 

The proof of \cite[Theorem 4.3]{FW} is incorrect, or at least incomplete;  
the  Hopf algebras $u_q^+({\mathfrak {g}})$ 
are not in general quasitriangular.
However,  Hochschild support variety theory \cite{BS,FW2}
provides an alternative proof
of this statement: 
By \cite[Lemma~6.3]{FW2} or \cite[\S2]{LZ}, 
the representation type of $u_q^+({\mathfrak{g}})$ is the same as
that of $u_q^{>0}({\mathfrak{g}})$, and we may apply 
\cite[Theorem~5.4]{FW2} to the latter algebra. 
The remainder of the proofs of the statements in \cite[Section 4]{FW} 
are correct, 
as the Hopf algebras $u_q({\mathfrak{g}})$ are quasitriangular.

We note that our results in \cite{FW2} almost exclusively use the
Hochschild support variety theory, and do not require any additional hypotheses. 
This theory is more general, and does not take
advantage of the tensor products of modules one has at hand for a Hopf algebra. 
One may obtain stronger results
by taking advantage of tensor products of modules,
yet the constructions are necessarily one-sided in general.
The one-sidedness affects results 
if the tensor product is not commutative up to isomorphism. 



\end{document}